\documentclass[12pt]{article}
\usepackage{fullpage}
\usepackage{amstex}
\usepackage{amssymb}
\usepackage{theorem}
\newtheorem{theorem}{Theorem}[section]
\newtheorem{lemma}[theorem]{Lemma}
{\theorembodyfont{\rmfamily}
\newtheorem{definition}[theorem]{Definition}%[section]
\newtheorem{remark}[theorem]{Remark}%[section]
}

\def\dom{{\rm dom}}
\def\l{{\langle}}
\def\r{{\rangle}}
\def\del{\delta}
\def\ome{\omega}
\def\mtp{\mathcal P}
\def\upr{\upharpoonright}
\def\bks{\backslash}	
	
\def\kap{\kappa}	
\def\lam{\lambda}

\def\ome{\omega}
\def\nek{,\ldots,}
\def\sig{\sigma}
\def\alp{\alpha}
\def\bet{\beta}
\def\orho{\overline\rho}

\title{Wide Gaps with Short Extenders}
\author{Moti Gitik\\
School of Mathematical Sciences\\
Tel Aviv University\\
 Tel Aviv 69978, Israel\\
gitik@@ math.tau.ac.il}
\date{}
\begin{document}
\renewcommand{\thefootnote}{$*$}

\maketitle
\begin{abstract}
The paper is a continuation of [Gi].
Extending the methods of [Gi] we show the
following:
Let $\kappa = \bigcup_{n<\omega}
\kappa_n$, $\kappa_0<\kappa_1<\cdots <
\kappa_n < \cdots$\\
 (1) if each $\kappa_n$
carries an extender of the length of the
first Mahlo  above $\kappa_n$, {\it
then} for every $\lambda > \kappa$ there is
a generic extension satisfying $2^\kappa
\ge \lambda$;\\
(2) if each $\kappa_n$ carries an extender
of the length of the first fixed point of
the $\aleph$-function above $\kappa_n$ of
order $n$, {\it then} for every $\lambda$,
$\kappa < \lambda<$ least inaccessible above
$\kappa$ there is a generic extension
satisfying $2^\kappa \ge \lambda$.
\end{abstract}

\baselineskip=18pt
\setcounter{section}{-1}
\renewcommand{\labelenumi}{(\theenumi)}
\renewcommand{\labelenumiii}{(\theenumiii)}
\section{Introduction}
We would like to extend the method of
[Gi]  to deal with arbitrary gaps
between a singular cardinal $\kappa$ and
its power.  The main structure of [Gi] was
as follows:\\
Let $\kappa = \bigcup_{n<\omega} \kappa_n$,
$\kappa_0 < \kappa_1 < \cdots < \kappa_n <
\cdots$ and we want to make $2^\kappa \ge
\kappa^{+\delta + 1}$ for some $\delta <
\kappa_0$.
\begin{enumerate}
\item For every $\nu \le \delta$ we use $\langle
\kappa_n^{+n+\nu + 1} \mid n < \omega\rangle$ as
a sequence below $\kappa$ responsible for
$\kappa^{+\nu + 1}$.
\item Sequences of elementary submodels of
some $H(\chi)$ ($\chi$ big enough) $\langle
A^{0\tau} \mid \tau \le \delta \rangle$
with $A^{0\tau} \subseteq A^{0\tau'}$
$(\tau \le \tau' \le \delta)$ and
$|A^{0\tau}| = \kappa^{+\tau+1}$ were used
to gradually  shrink the number of possible
choices of generic $\omega$-sequences.
Actually, the first submodels $A^{00}$'s
are the most important, since they have the
size $\kappa^+$ which insures eventually
the $\kappa^{++}$-c.c. of the final
forcing.  It was essential that $A^{0\tau}$'s
include $\delta$ and we can deal at once
with $\langle \sup (A^{0\tau} 
\cap \kappa^{+\nu + 1} )\mid \nu \le
\delta\rangle$.
\end{enumerate}

Stretching the forcing of [Gi] slightly, we can
deal with $\delta$'s above $\kappa$ but
below $\kappa^{++}$ and to make $2^\kappa$
as big as we like below $\kappa^{+
(\kappa^{++})} = \aleph_{\kappa +
\kappa^{++}}= \aleph_{\kappa^{++}}$.  Thus for example if $\delta
= \kappa^+ + 1$, then in (1) we may use $\langle
\kappa_n^{+\kappa_n^+ + \kappa_n^+ + 1}\mid
n<\omega\rangle$ instead of $\langle
\kappa_n^{+n+\delta + 1} \mid n < \omega\rangle$
and take care of $\kappa^{+\kappa^+ + 1}$.
This will leave enough room for $\nu$'s below
$\delta$ as well.  These is no problem with
(2), since the number of cardinals we are
dealing with is still below $\kappa^{++}$.
In Section 2 we will do it under weaker
assumptions.

But once we are above $\kappa^{++}$, (1)
and especially (2), require essential
changes.
  
We show in the first section the following:

%\begin{theorem}
\medskip
\noindent
{\bf Theorem 1}\footnote{S. Shelah
suggested to push this down combining the
forcing of the theorem with [Gi-Ma, Sec. 2]
or [Gi1, Sec. 4]. It is possible this way
to obtain models with no Mahlo below
$\kappa$ and $2^\kappa$ arbitrary large.}
\quad
{\it Let $\kappa = \cup_{n<\omega}\kappa_n$,
$\kappa_0 < \kappa_1 < \cdots < \kappa_n <
\cdots$ and each $\kappa_n$ carries an
extender of the length of the first
Mahlo above $\kappa_n$. \underline{Then}
for every $\lambda > \kappa$ there is a
cofinality preserving extension, with
the exception of the successors of
inaccessibles above $\kappa$, not adding
new bounded subsets to $\kappa$ and
satisfying $2^\kappa \ge \lambda$.}
%\end{theorem}

In the third section we show the following:
%\begin{theorem}
\medskip

\noindent
{\bf Theorem 2}\quad {\it Let $\kappa=\bigcup_{n<\omega} \kappa_n$,
$\kappa_0 < \kappa_1 < \cdots < \kappa_n <
\cdots$ and each $\kappa_n$ carries an
extender of the length of the first fixed
point of the $\aleph$-function above
$\kappa_n$ of order $n$. \underline{Then} for
every $\lambda>\kappa$ and below the first
inaccessible above $\kappa$ there is a
cofinality preserving extension not adding
new bounded subsets to $\kappa$ and
satisfying $2^\kappa \ge \lambda$.}
%\end{theorem} 

\section{Getting an arbitrary gap}
Our purpose here will be to define the
forcing to prove Theorem 1.  We will use
forcings similar to those of [Gi, Sec.~4]
with changes for overcoming (1) and (2)
above.
 
Fix an ordinal $\delta > 1$.

\begin{definition}
The set ${\cal P}'$ consists of pairs
$\langle\langle A^{0\tau}, A^{1\tau}
\rangle \mid \tau \le \delta\rangle$ so
that the following holds:
\begin{enumerate}
\item for every $\tau \le \delta$
$A^{0\tau}$ is an elementary submodel of
$\langle H(\kappa^{+\delta + 2}),
\epsilon, \langle \kappa^{+i} \mid i \le
\delta + 2\rangle\rangle$ such that
\begin{enumerate}
\item $|A^{0\tau}| = \kappa ^{+\tau + 1}$
and $A^{0\tau} \supseteq \kappa^{+\tau +
1}$ unless for some $n < \omega$ and an
inaccessible $\tau'$, $\tau = \tau' + n$
and then $|A^{0\tau}| = \kappa^{+\tau}$ and
$A^{0\tau} \supseteq \kappa^{+\tau}$
\item ${}^{|A^{0\tau}| >}
\!\!A^{0\tau} \subseteq A^{0\tau}$
\end{enumerate}
\item for every $\tau < \tau' \le \delta$,
 $A^{0\tau} \subseteq A^{0\tau'}$
\item for every $\tau \le \delta$,
 $A^{1\tau}$ is a set of at most
$\kappa^{+\tau + 1}$ elementary submodels
of $A^{0\tau}$ so that
\begin{enumerate}
\item $A^{0\tau} \in A^{1\tau}$
\item if $B,C \in A^{1\tau}$ and $B
\subsetneqq C$ then $B \in C$
\item if $B \in A^{1\tau}$ is a successor
point of $A^{1\tau}$ then $B$ has at most
two immediate predecessors under the
inclusion and is closed under
$\kappa^{+\tau}$-sequences.
\item let $B \in A^{1\tau}$ then either $B$
is a successor point of $A^{1\tau}$ or $B$
is a limit element and then there is a
closed chain of elements of $B \cap
A^{1\tau}$ unbounded in $B \cap A^{1\tau}$
and with limit $B$.
\item for every $\tau',\tau \le \tau' \le
\delta$, $A \in A^{1\tau}$ and $B \in
A^{1\tau'}$ either $B \supseteq A$ or there
are $\ell < \omega$ and $\tau'_1,\tau'_2,\ldots
, \tau'_\ell$, $\tau \le \tau'_1 \le \cdots
\le \tau'_\ell \le \delta$, $B_1 \in A \cap
A^{1\tau'_1}, \ldots, B_\ell \in A \cap
A^{1\tau'_\ell}$ such that
$$B \cap A = B_1 \cap \cdots \cap B_\ell
\cap A\ ,$$
if $\tau = \tau'$, then we can pick
$\tau'_1$ (and hence all the rest) above
$\tau$.

\item let $A$ be an elementary submodel of
$H(\kappa^{+\delta + 2})$ of cardinality
$|A^{0\tau}|$, closed under
$<|A^{0\tau}|$-sequences,
$|A^{0\tau}| \in A$ and including
$\langle\langle A^{0\tau'}, A^{1\tau'}
\rangle \mid \tau' \le \delta \rangle$ as
an element, for some $\tau \le \delta$.
Then for every $\tau'$, $\tau\le \tau' \le
\delta$ and $B \in A^{1\tau'}$ either $B
\supseteq A$ or there are $\tau'_1,\ldots ,
\tau'_\ell$, $\tau \le \tau'_1 \le \cdots \le
\tau'_\ell \le \delta$, $B_1 \in A \cap
A^{1\tau'_1}, \ldots, B_\ell \in A \cap
A^{1\tau'_\ell}$ such that
$$B \cap A = B_1 \cap \cdots \cap B_\ell
\cap A\ .$$
\end{enumerate}
\end{enumerate}
\end{definition}

The present definition of ${\cal P}'$ differs from
those of [Gi, 4.1]
by two additional conditions (f) and (g).
They are desired in order to overcome the
difficulty (2).  This way the number of
possible intersections (and actually
intersections themselves) is controlled.

The addition of (f) and (g) makes the proof
of distributivity of the forcing more
involved and  we shall further concentrate
on this matter.

The definition of order 1.2, 1.3 and 1.4
repeats the corresponding ones in [Gi,
Sec.~4].

Let for $\tau\le \del$  $A^{1\tau}_{in}$  be the
set $\{B\cap B_1\cap \cdots\cap B_n\mid B\in
A^{1\tau}$, $n<\ome$  and $B_i\in A^{1\rho_i}$ for
some $\rho_i,\tau <\rho_i\le\del$ for
every $i,1\le i\le n\}$. 

\begin{definition}
  Let
$\langle\langle A^{0\tau},A^{1\tau}\rangle\mid\tau\le\del
\rangle$  and $\langle\langle B^{0\tau},
B^{1\tau}\rangle\mid \tau\le\del\rangle$  be elements of
$\mtp'$.  Then $\langle\langle A^{0\tau},A^{1\tau}\rangle
\mid\tau\le\del\rangle\ge\langle\langle
B^{0\tau},B^{1\tau}\rangle\mid\tau\le\del\rangle$
iff for every $\tau\le\del$
\begin{enumerate}
\item $A^{1\tau}\supseteq B^{1\tau}$
\item for every $A\in A^{1\tau}$ either
\begin{enumerate}
\item $A \supseteq B^{0\tau}$ or
\item $A \subset B^{0\tau}$ and then $A \in
B^{1\tau}$ or
\item $A \not\supseteq B^{0\tau}$,
$B^{0\tau} \not\supseteq A$ and then $A \in
B_{in}^{1\tau}$.
\end{enumerate}
\end{enumerate}

\end{definition}

\begin{definition}
  Let
$\tau\le\del$.  Set $\mtp'_{\ge\tau}=\{\langle\langle
A^{0\rho},A^{1\rho}\rangle\mid\tau\le\rho\le\del\rangle
\mid\exists\langle\langle
A^{0\nu},A^{1\nu}\rangle\mid\nu <\tau\rangle$
$\langle\langle A^{0\nu},A^{1\nu}\rangle\mid\nu
<\tau\rangle^\frown\langle\langle A^{0\rho},
A^{1\rho}\rangle\mid\tau\le\rho\le\del\rangle\in\mtp'\}$.

Let $G(\mtp'_{\ge\tau})\subseteq\mtp'_{\ge\tau}$
be generic.  Define $\mtp'_{<\tau}=\{\langle\langle
A^{0\nu},A^{1\nu}\rangle \mid\nu <\tau\rangle\mid\exists
\langle\langle A^{0\rho},A^{1\rho}\rangle\mid\tau\le\rho
\le\del\rangle\in G(\mtp'_{\ge\tau})$ $\langle\langle
A^{0\nu},A^{1\nu}\rangle\mid\nu <\tau\rangle\raise
3pt\hbox{$\frown$}\langle
\langle A^{0\rho},A^{1\rho}\rangle\mid\tau\le\rho\le\del
\rangle\in\mtp'\}$.
\end{definition}

The following lemma is obvious

\begin{lemma}
 $\mtp'\simeq\mtp'_{\ge\tau}*
\underset{\sim}{\mtp'}_{\raise .5em\hbox{$\scriptstyle
<\tau$}}\quad (\tau\le\del)$.
\end{lemma}

Let us now define the main preparation
forcing ${\cal P}$.  The definition repeats
mainly the corresponding definition of
[Gi], 4.14. The only changes are designed
to overcome (1), since we are now probably dealing
with a large number of cardinals
above $\kappa$.

Thus we like cardinals to correspond to
cardinals, regular cardinals to regular cardinals, and limit cardinals
to correspond to limit cardinals.  This puts
some limitations.  In particular,
inaccessibles should correspond to
inaccessibles, inaccessibles of order 1
(i.e. inaccessibles which are limits of
inaccessibles) should correspond to
inaccessibles of order $\ge 1$, etc.  We
shall arrange in a moment sets of
inaccessibles below $\kappa$ corresponding
to those above $\kappa$.
But first we destroy all Mahlo cardinals
above $\kappa$ by forcing for each such
cardinal a  club avoiding inaccessibles.

Fix $n<\omega$.  We now want to define 
``good" inaccessibles.  This notion will be
similar to the notion of good ordinals used
in [Gi,2.8]. Let $\chi_n$ denote the least
Mahlo cardinal above $\kappa_n$.  For every
$k \le n$ we consider the structure $\mathfrak
{a}_{n,k} = \l H(\chi_n^{+ k + 2} ), \epsilon,
E_n, 0,1,\ldots, \alpha, \ldots \mid \alpha
< \kappa_n^{+k+2}\r$. For an ordinal $\xi <
\chi_n$ let $tp_{n,k} (\xi)$ be the type
realized by $\xi$ in $\mathfrak{a}_{n,k}$.

The following lemma is obvious:
\medskip

\noindent
{\bf Lemma 1.4.1} {\it There are stationary
many below $\chi_n$ inaccessible cardinals
$\xi$ so that
\begin{itemize}
\item[{\rm (1)}] the set $\{\delta < \xi
\mid tp_{n,n}(\delta) = tp_{n,n} (\xi) \}$
is unbounded in $\xi$,
\item[{\rm (2)}] $\xi = A\cap \chi_n$ for
some $A \prec \mathfrak{a}_{n,n}$.
\end{itemize}
}
\smallskip

We fix an inaccessible $\xi_n$ satisfying
the conclusion of the lemma.

\medskip
\noindent
{\bf Definition 1.4.2}  An inaccessible
cardinal $\delta < \xi_n$ is called
$k$-good (for some $k \le n$) if $tp_{n,k}
(\delta) = tp_{n,k} (\xi_n)$.
\medskip

The next lemma is obvious.
\medskip

\noindent
{\bf Lemma 1.4.3} {\it If $\delta$ is
$k$-good for some $k > 1$, then $\delta$ is
the limit of $k-1$-good inaccessibles.}
\smallskip

Further we shall use only the restriction
of the extender $E_n$ to $\xi_n$.  Also
each inaccessible $\lambda > \kappa$ will
correspond to a sequence $\langle \delta_n
\mid n < \omega \r$, where every $\delta_n$
is $k_n$-good inaccessible and $k_0 < k_1 <
\cdots < k_n < \cdots$.

\begin{definition}
The set ${\cal P}$ consists of sequences of
triples $\langle\langle A^{0\tau},
A^{1\tau}, F^\tau \r \mid \tau\le \delta
 \rangle$ so that the following holds:

\begin{enumerate}
\setcounter{enumi}{-1}
\item $\l\l A^{0\tau}, A^{1\tau} \r\mid \tau
\le \delta \r\in {\cal P}'$
\item for every $\tau_1 \le \tau_2\le
\delta$, $F^{\tau_1} \subseteq F^{\tau_2}
\subseteq {\cal P}^*$
\item for every $\tau \le \delta$, $F^\tau$
is as follows:
\begin{enumerate}
\item $|F^\tau| = |A^{0\tau}|$
\item for every $p = \l p_n \mid n < \omega
\r \in F^\tau$ if $n < \ell (p)$ then every
$\alpha$ appearing in $p_n$ is in $A^{0\tau}
\cup \{|A^{0\tau}|\}$; if $n \ge \ell (p)$
and $p_n = \l a_n, A_n, f_n\r$ then every
$\alpha$ appearing in $f_n$ is in
$A^{0\tau} \cup \{ |A^{0\tau}|\}$
 and
\begin{enumerate}
\item $\dom a_n \cap On \subseteq
(A^{0\tau} \cap \kappa^{+ \delta + 2}) \cup
\{ |A^{0\tau}| \}$
\item $\dom a_n \backslash On$ consists of
elements of the following sets: $\{ B
\subseteq A^{0\tau}\mid \kappa^+ \le |B| <
|A^{0\tau}|\}$, $A^{1\tau}$ and
$A_{in}^{1\tau}$ such that the elements of
the last two sets are closed under $>
|A^{0\tau}| $-sequences of
its elements.  If $\tau = 0$, then the
first set is empty.
\end{enumerate}
\item  the largest (under
inclusion) element of $\dom a_n \backslash
On$ belongs to $A^{1\tau}$ and every
element of $\dom a_n$ belongs to it.

Let us further denote this element
as $\max^1 (p_n)$ or $\max^1 (a_n)$.
\item if $B \in \dom a_n \backslash On$,
then $a_n(B)$ is an elementary submodel of
$a_{n,k_n}$ of Section 2  of [Gi] with $3 \le k_n
\le n$, including also $\delta$ as a
constant.  We require that 
\item[(d1)] if $|B|$ is
a successor cardinal, then $|a_n (B)| =
\kappa_n^{+\tau' + 1}$ and
${}^{\kappa_n^{+\tau'}} (a_n(B)) \subseteq
a_n(B)$, where $\tau'$, $\kappa_n^+ < \tau'
<$ $\xi_n$
is $k_n$-good cardinal, $k_n$'s are
increasing with $n$ and $a_n (|B|) =
\kappa_n^{+\tau' + 1}$.
\item[(d2)] if $|B|$ is an inaccessible
cardinal, then $|a_n(B)|$ is a $k_n$-good
inaccessible with $k_n$'s increasing with
$n$, ${}^{|a_n(B)|>} a_n (B) \subseteq
a_n(B)$ and $a_n (|B|) = |a_n (B)|$.
\item if $B \in \dom a_n \backslash On$
and $\alpha \in \dom a_n \cap A^{0\tau}$
then $a_n(\alpha) \in a_n(B)$ iff $\alpha
\in B$
\item if $B,C\in \dom a_n\backslash On$
then 
\item[(f1)] $B \in C$ iff $a_n (B) 
\in a_n (C)$
\item[(f2)] $B \subset C$ iff $a_n (B)
\subset a_n (C)$.
\\
The next condition deals with cofinalities
correspondence
\item (i) if $\alpha \in \dom a_n$ and
$cf\alpha \le \kappa^+$ then $cf a_n
(\alpha) \le \kappa_n^{+n + 1}$.
\begin{itemize}
\item[(ii)] if $\alpha \in \dom a_n$
and $cf\alpha = \kappa^{+\rho}$ for some $\rho$,
$1 \le \rho \le \delta + 1$ then
$\kappa^{+\rho} \in \dom a_n$, $cf a_n
(\alpha) = a_n(\kappa^{+\rho})$ and for
every $B \in \dom a_n \backslash On$ of
cardinality $\kappa^{+\rho}$\quad  $|a_n (B)|=a_n
(\kappa^{+\rho})$.
\item[(iii)] if $\alpha \in \dom a_n$
is an inaccessible, then $a_n (\alpha)$ is
$k_n$-good inaccessible, with $k_n$'s
increasing with $n$.
\end{itemize}
%
%
%***************************************

%INSERT PAGES 15-17
%
%***************************************
%
%\begin{itemize}
\item[(h)] if $p\in F^\tau$ and $q\in\mtp^*$  is
equivalent to $p$ $\quad(q\leftrightarrow p)$  with
witnessing sequence $\langle k_n\mid
n<\ome\rangle$  starting with $k_0\ge 4$  then
$q\in F^\tau$.
\item[(i)] if $p=\langle p_n\mid
n<\ome\rangle\in F^\tau$  and $q=\langle q_n\mid
n<\ome\rangle\in\mtp^*$  are such that
\begin{itemize}
\item[(i)] $\ell(p)=\ell (q)$ 
\item[(ii)] for every $n<\ell (p)$ $p_n=q_n$ 
\item[(iii)] for every $n\ge\ell (p)$ $a_n=b_n$
and $\dom g_n\subseteq A^{0\tau}$ where $p_n=\langle
a_n,A_n,f_n\rangle$, $q_n=\langle b_n,B_n,g_n\rangle$
\end{itemize}
\item[] then $q\in F^\tau$.
%\end{itemize}
%\begin{itemize}
\item[(k)] if $p=\langle p_n\mid n<\ome\rangle\in F^\tau$
$q=\langle q_n\mid n<\ome\rangle\in\mtp^*$  are such that 
%\end{itemize}
\begin{itemize}
\item[(i)] $\ell (q)\ge\ell (p)$
\item[(ii)] for every $n\ge\ell (q)$ $p_n=q_n$
\item[(iii)] every $\alp$  appearing in $q_n$
for $n<\ell (q)$  is in $A^{0\tau}$
\item[] then $q\in F^\tau$.
\end{itemize}

The meaning of the last two conditions is that
we are free to change inside $A^{0\tau}$ all the
components of $p$ except $a_n$'s.
%\begin{itemize}
\item[(l)] for every $q\in F^\tau$  and $\alp
\in A^{0\tau}$  there is $p\in F^\tau$
$p=\langle p_n\mid n<\ome\rangle$, $p_n=\langle
a_n,A_n,f_n\rangle$  $(n\ge\ell (p))$ such that
$p\phantom{!}^*\!\!\geq q$  and $\alp\in\dom a_n$  starting
with some $n_0<\ome$. 
\item[(m)] for every $q\in F^\tau$  and $B\in
A^{1\tau}\cup A^{1\tau}_{in}$  as in (b)(ii),
there is $p\in F^\tau$ $p=\langle p_n\mid
n<\ome\rangle$,  $p_n=\langle a_n,A_n,f_n\rangle$  
$(n\ge\ell (p))$  such that $p\ge^* q$  and
$B\in\dom a_n$  starting with some $n_0<\ome$.
Also, this $p$ is obtained from $q$ by adding only
$B$ and the ordinals needed to be added after
adding $B$. 
\item[(n)] Let $p,q\in F^\tau$  be so that 
\begin{itemize}
\item[(i)] $\ell (p)=\ell (q)$
\item[(ii)] $\max^1(p_n)=\max^1(p_n)$,
$\max^1(q_n)=\max^1(q_m)$ and $\max^1(q_n)\in\dom
a_n$, where $n,m\ge\ell (p)$, $p_n=\langle
a_n,A_n,f_n\rangle$, $q_n=\langle b_n,B_n,g_n\rangle$ 
\item[(iii)] $p_n=q_n$ for every $n<\ell (p)$
\item[(iv)] $f_n,g_n$  are compatible for every
$n\ge\ell (p)$
\item[(v)] $a_n\upr\max^1(q_n)\subseteq b_n$
for every $n\ge\ell(p)$, where
$$a_n\upr B=\{\langle t\cap B,s\cap
a_n(B)\rangle \mid \langle t,s\rangle\in a_n\}$$
then the union of $p$ and $q$ is in $F^\tau$ where
the union is defined in obvious fashion taking
$p_n\cup q_n$  for $n<\ell (p)$, we take at each
$n\ge\ell (p)$ $a_n\cup b_n$, $f_n\cup g_n$  etc. 
\end{itemize}
%\end{itemize}
%\begin{itemize}
\item[(o)] there is $F^{\tau *}\subseteq F^\tau$
dense in $F^\tau$  under $\le^*$  such that
every $\le^*$-increasing sequence of elements of
$F^{\tau *}$  having the union in $\mtp^*$  has
it also in $F^\tau$.
We require that $F^{\tau*}$  will be closed
under the equivalence relation
$\leftrightarrow$.
\item[(p)] let $p=\langle p_n\mid n<\ome\rangle
\in F^\tau$ and $p_n=\langle a_n,A_n,f_n\rangle$
$(\ell (p)\le n<\ome)$.  If for every $n$,
$\ome >n\ge\ell (p)$ $B\in\dom a_n\bks On$,
$|B|=\kap^{+\tau +1}$ or  $B\in A^{1\tau'}$ for some
$\tau'\le\tau$, \underline{then}
$p\upr B\in F^{\tau'}$, where $p\upr
B=\langle p_n\upr B\mid n<\ome\rangle$ and
for every $n<\ell (p)$ $p_n\upr B$  is the
usual restriction of the function $p_n$  to $B$;
if $n\ge\ell (p)$ then $p_n\upr B=\langle a_n\upr 
B, B_n,f_n\upr B\rangle$ with $a_n\upr
B$  defined in (n)(v), $f_n\upr B$  is the
usual restriction and $B_n$  is the projection
of $A_n$  by $\pi_{\max p_n,B}$.
\item[(q)] let $p=\langle p_n\mid n<\ome\rangle\in
F^\tau$, $p_n=\langle a_n,A_n,f_n\rangle$ and
$A^{0\tau}\not\in \dom a_n$ $(\ome >n\ge\ell (p))$.
Let $\langle \sig_n\mid\ome >n\ge\ell
(p)\rangle$ be so that 
\begin{itemize}
\item[(i)] $\sig_n\prec\mathfrak{a}_{n,k_n}$ 
and $|\sigma_n|$ is $k_n$-good for
every $n\ge\ell (p)$
\item[(ii)] $\langle k_n\mid n\ge\ell (p)\rangle$
is increasing
\item[(iii)] $k_0\ge 5$
\item[(iv)] ${}^{|\sigma_n| >}\!\sig_n\subseteq
\sig_n$  for every
$n\ge\ell (p)$ 
\item[(v)] $rnga_n\in\sig_n$  for every
$n\ge\ell(p)$. 

Then the condition obtained from
$p$  by adding $\langle A^{0\tau},\sig_n\rangle$
to each $p_n$  with $n\ge\ell(p)$  belongs to
$F^\tau$.
\end{itemize}
\item[(r)] if $A$ is an elementary submodel of
$H(\kap^{+\del +2})$  of a regular cardinality
$\kap^{+\rho}$, closed under $<\kap^{+\rho}$-sequences
and including $\langle\langle A^{0\tau'},
A^{1\tau'}\rangle \mid\tau'\le\del\rangle$ 
for some $\rho <\tau$, then $A$ is addable to
any $p\in F^\tau\cap A$, with the maximal element
of $\dom a_n$'s $A^{0\tau}$,  i.e. $A\cap
A^{0\tau}$ can be added to $p$  remaining in
$F^\tau$. 
%\end{itemize}
\end{enumerate}
\end{enumerate}
\end{definition}

\begin{definition}  Let
$\langle\langle A^{0\tau}$, $A^{1\tau},
F^\tau\rangle\mid\tau\le\del\rangle$  and
$\langle\langle B^{0\tau},B^{1\tau},G^\tau\rangle\mid
\tau\le\del\rangle$ be in $\mtp$.  We define
$$\langle\langle A^{0\tau}, A^{1\tau},F^\tau\rangle\mid
\tau\le\del\rangle >\langle\langle B^{0\tau},
B^{1\tau},G^\tau\rangle\mid\tau\le\del\rangle$$
iff
\begin{itemize}
\item[(1)] $\langle\langle
A^{0\tau},A^{1\tau}\rangle\mid\tau\le\del\rangle
>\langle\langle B^{0\tau},B^{1\tau}\rangle
\mid\tau\le\del\rangle$ in $\mtp'$
\item[(2)] for every $\tau\le\del$
\end{itemize}
\begin{itemize}
\item[(a)] $F^\tau\supseteq G^\tau$
\item[(b)] for every $p\in F^\tau$  and $B\in
B^{1\tau}\cup B^{1\tau}_{in}$ if for every
$n\ge\ell(p)$ $B\in\dom a_n$  then $p\upr
B\in G^\tau$, where the restriction is as defined
in 1.5(p), $p=\langle p_n\mid n<\ome\rangle$,   
$p_n=\langle a_n, A_n,f_n\rangle$  for $n\ge
\ell (p)$.
\end{itemize}
\end{definition}

\begin{definition}
  Let
$\tau\le\del$.  Set $\mtp_{\ge\tau}=\{\langle
A^{0\rho}, A^{1\rho}, F^\rho\rangle \mid
\tau\le\rho\le\del\rangle\mid\exists\langle\langle
A^{0\nu}, A^{1\nu}, F^\nu\rangle\mid\nu
<\tau\rangle$ $\langle\langle A^{0\nu},A^{1\nu},
F^\nu\rangle\mid\nu <\tau\rangle\frown\langle\langle
A^{0\rho}, A^{1\rho},F^\rho\rangle\mid\tau\le\rho\le
\del\rangle\in\mtp\}$.
\medskip

Let $G(\mtp_{\ge\tau})\subseteq\mtp_{\ge\tau}$
be generic.  Define $\mtp_{<\tau}=\{\langle\langle
A^{0\nu},A^{1\nu},F^\nu\rangle\mid\nu
<\tau\rangle,\mid\exists\langle\langle A^{0\rho}$,
$A^{1\rho}, F^\rho\rangle \mid\tau\le\rho\le\del\rangle\in
G(\mtp_{\ge\tau})$ $\langle\langle A^{0\nu}$,
$A^{1\nu}, F^\nu\rangle\mid\nu
<\tau\rangle\frown\langle\langle
A^{0\rho},A^{1\rho}$, $F^\rho\rangle\mid\tau\le\rho\le\del
\rangle\in\mtp\}$.
\end{definition}

The following lemma is obvious

\begin{lemma}
 $\mtp\simeq\mtp_{\ge\tau}
*\underset{\sim}{\mtp}_{\raise .5em\hbox{$\scriptstyle
<\tau$}}$ for every $\tau\le\del$.
\end{lemma}

\begin{lemma}
For every $\tau \le \delta$, ${\cal P}_{\ge
\tau}$ is $\kappa^{+\tau +
2}$-strategically closed.  Moreover, if
there are an inaccessible $\tau'\le \tau$
and $n < \omega$ such that $\tau = \tau' +
n$, then ${\cal P}_{ \ge \tau}$ is
$\kappa^{+\tau + 1}$-strategically closed.
\end{lemma}

\medskip
\noindent
{\bf Proof.} The proof is similar to
[Gi,4.18] and we concentrate on new points
caused by the additions to the definition
of ${\cal P}'$ here, i.e. (f) and (g) of
1.1.

Fix $\tau \le \delta$.  Let $\l\l
A_i^{0\rho}, A_i^{1\rho}, F_i^\rho\r \mid
i<i^*\r$ be an increasing sequence of
conditions in ${\cal P}_{\ge \tau}$ already
generated by playing the game.  We need
to define the move $\l\l A_{i^*}^{0\rho},
A_{i^*}^{1\rho}, F_{i^*}^\rho\mid \tau \le
\rho \le \delta\r$ of Player I at stage
$i^*$.

Define this triple by induction on $\rho$.
The definition of $F_{i^*}^\rho$ completely
repeats
the one in [Gi, 4.18].  So we
deal only with $\l A_{i^*}^{0\rho},
A_{i^*}^{1\rho}\r$. 
\medskip

\noindent
{\bf Case 1}\quad $i^*$ is a limit and $cf i^* =
\kappa^{+\tau + 1}$ (or $\kappa^{+\tau}$, if
$\tau = \tau' + n$ for inaccessible $\tau'
\le \tau$ and $n < \omega$).

 We set
$A_{i^*}^{0\tau} = \bigcup_{i < i^*}
A^{0\tau}_i$ and $A_{i^*}^{1\tau} =
\bigcup_{i < i^*} A_i^{1\tau}
\bigcup \{ A_{i^*}^{0\tau}\}$.  Let $\rho
\in (\tau,\delta]$.  Set $\widetilde
A_{i^*}^{0\rho}$ to be the closure under
the Skolem functions and
$\kappa^{+\rho}$-sequences (or $<
\kappa^{+\rho}$-sequences, if $\rho =
\rho^* + n$ for an inaccessible $\rho^*$
and $n < \omega$) of $\l\l A_i^{j \rho'}
\mid i < i^*\r \mid\tau \le \rho' \le
\delta\r$ $(j \in 2)$ and $\l
A_{i^*}^{1\rho'} \mid \tau \le \rho'
<\rho\r$.  Define $A_{i^*0}^{0\rho}$ to be
the closure under the Skolem functions
and $\kappa^{+\rho}$-sequences (or $<
\kappa^{+\rho}$-sequences, if $\rho^*=\rho
+ n$ for an inaccessible $\rho^*$ and $n <
\omega$) of $\widetilde A_{i^*}^{0\rho}$,
$\l F^{\rho'}_i \mid \rho \le \rho' \le
\delta, i < i^*\r$ and $\l F_i^{\rho' *}
\mid \rho \le \rho' \le \delta, i < i^*\r$.
Let $\kappa (\rho) = \kappa^{+\rho}$ if
there is an inaccessible $\rho^* \le \rho$
and $n < \omega$ such that $\rho = \rho^* +
n$, and let $\kappa (\rho) = \kappa^{+\rho
+ 1}$ otherwise. For a limit $\alpha$, $0 <
\alpha < \kappa (\rho)$ let $A_{i^*\alpha}
^{0\rho} = \bigcup_{\alpha' <
\alpha} A^{0\rho}_{i^* \alpha'}$.  Let
$A_{i^* \alpha + 1}^{0\rho}$ be the closure
of $A^{0\rho} _{i^* \alpha} \cup \{ A_{i^*
\alpha}^{0\rho}\}$ under the Skolem
functions and $< \kappa(\rho)$-sequences,
for every $\alpha < \kappa (\rho)$.  We set
$A^{0\rho}_{i^* \alpha} =
\bigcup_{\alpha < \kappa (\rho)}
A_{i^*\alpha}^{0\rho}$ and $A_{i^*}^{1\rho}
= \bigcup_{i < i^*} A_i^{1\rho} \cup
\{ A^{0\rho}_{i^* \alpha} \mid \alpha <
\kappa (\rho)\}\cup \{ A_{i^*}^{0\rho}\}$.

\medskip
\noindent
{\bf Case 2}\quad $i^*$ is a successor ordinal
or $i^*$ is a limit ordinal of cofinality
$ \kappa^{+\tau + 1}$ or less than
$\kappa^{+\tau}$ if $\tau = \tau' + n$ for
an inaccessible $\tau'$ and $n < \omega$.

In this case we treat $\tau$ in the same
way as any other $\rho \in (\tau,\delta]$
in the previous case.  The definition for
$\rho \in (\tau,\delta]$ is as
in Case 1.

Let us show now that such defined $\l\l
A_{i^*}^{0\rho}, A_{i^*}^{1\rho}\r \mid
\tau \le \rho \le \delta\r$ is in ${\cal
P}'$.  Basically, we need to check the
conditions (f) and (g) of Definition 1.1.

We start with (f). Let $\tau \le \rho \le
\rho' \le \delta$, $A \in A_{i^*}^{1\rho}$
and $B \in A_{i^*}^{1\rho'}$. If $A \in
A_i^{1\rho}$ and $B \in A_{i'}^{1\rho'}$
for some $i$, $i' < i^*$, then we use (f)
for $\l\l A_{\overline i}^{0\nu},
A_{\overline i}^{1\nu}\r \mid \tau \le
\nu\le \delta\r$ where $\overline i = \max
(i,i')$.  It provides $\rho \le \tau'_1\le
\cdots \le \tau'_\ell \le \delta$, $B_1
\in A \cap A_{\overline i}
^{1\tau'_1},\ldots, B_\ell \in A \cap
A_{\overline i} ^{1 \tau'_\ell}$ such that
$B \cap A = B_1 \cap \cdots \cap B_\ell
\cap A$.  Now, since $A_{\overline i}^{1
\tau'_k} \subseteq A_{i^*}^{1 \tau'_k}$ for
every $1 \le k \le \ell$ we are done.

If $A \in A_i^{1\rho}$ for some $i < i^*$
and $B \in A_{i^*}^{1\rho'} \backslash
\bigcup_{i < i^*} A_{i'}^{1\rho'}$
then $B \supseteq \bigcup_{i' < i^*}
A_{i'}^{0 \rho'}$.  In particular, $B
\supseteq A_i^{0\rho'} \supseteq
A_i^{0\rho}$.  If $A \in A_{i^*}^{1\rho}
\backslash \bigcup_{i < i^*}
A_i^{1\rho}$ and $B \in A_{i'}^{1\rho'}$
for some $i' < i^*$, then we can use 1.1(g)
for  
$A,B$ and $\l\l A_{i'}^{0\tau'},
A_{i'}^{1\tau'} \r \mid \tau' \le \delta \r
\in {\cal P}'$. If $A \in A_{i^*}^{1\rho}
\backslash \bigcup_{i < i^*}
A_i^{1\rho}$ and $B \in A_{i^*}^{1\rho'}
\backslash \bigcup_{i < i^*} A_i^{1
\rho'}$, then either $B \supseteq A$ or $B
\subset A$ and in this case $\rho' = \rho$
and $B \in A$.

Now let us check the condition (g). Thus
let $A$ be an elementary submodel of
$H(\kappa ^{+\delta + 2})$ of cardinality
$|A^{0\rho}_{i^*}|$, closed under $<
|A_{i^*}^{0\rho}|$-sequences,
$|A_{i^*}^{0\rho}| \in A$ and including
$\l\l A_{i^*}^{0\tau'}, A_{i^*}^{1\tau'} \r
\mid \tau' \le \delta \r $ as an element,
for some $\rho \le \delta$.  Let $\tau' \in
[\rho,\delta]$ and $B \in
A_{i^*}^{1\tau'}$.  Suppose first that $B
\in A^{1\tau'}_{i'}$ for some $i' < i^*$.
Then, $\l\l A_{i'}^{0\nu}, A_{i'}^{1\nu} \r
\mid \nu \le \delta \r \in A$, since
$A_{i^*}^{0\tau} \subseteq A_{i^*}^{0\rho}
\subseteq A$ and the sequence $\l\l
A_{i'}^{0\nu}, A_{i'}^{1\nu} \r \mid \nu
\le \delta \r \in A_{i^*}^{0\tau}$.  So (g)
of 1.1 applies to $A,B$ and $\l\l
A_{i'}^{0\nu}, A_{i'}^{1\nu} \r \mid \nu
\le \delta \r$ and we are done.  Assume now
that $B \in A_{i^*}^{1\tau'} \backslash
\bigcup_{i < i^*} A_i^{1\tau'}$.  If $\tau'
= \rho$, then $B \in A$ since $A \supseteq
|A_{i^*}^{0\rho}| = \kappa (\rho)$,
$A_{i^*}^{0\rho} \in A$ and, so
$A_{i^*}^{0\rho} \subseteq A$. But either
$B = A_{i^*}^{0\rho}$ or $B 
\in A_{i^*}^{0\rho}$.  Suppose now that
$\tau'>\rho$.  If $\tau' \in A$, then
$A_{i^*}^{0\tau'} \in A$.  Recall that in
this case
$A_{i^*}^{0\tau'} = \bigcup_{\alpha
< \kappa(\tau')} A_{\alpha i^*}^{0\tau'}$
and $A_{i^*}^{1\tau'} = \bigcup_{i <
i^*} A_i^{1\tau'} \cup \{
A_{i^*\alpha}^{0\tau'} \mid \alpha < \kappa
(\tau') \} \cup \{A_{i^*}^{0\tau'}\}$. If
$B = A_{i^*}^{0\tau'}$ or $B =
A_{i^*\alpha} ^{0\tau'}$ for some $\alpha
\in A$ then we are done.  Suppose
otherwise. Then, let $B=
A_{i^*\alpha}^{0\tau'}$ for some $\alpha <
\kappa (\tau')$.  Set $\widetilde \alpha =
\min (A\backslash \alpha)$.  Then
$\widetilde \alpha \le \kappa (\tau')$ and
$A^{0\tau'}_{i^* \tilde\alpha} \in A$,
where $A_{i^*\kappa (\tau')}^{0\tau'} =
A^{0\tau'}_{i^*}$. But now $A \cap B = A
\cap 
A^{0\tau'}_{i^* \tilde\alpha}$, since the
chain $\l A_{i^*\beta}^{0\tau'} \mid \beta
< \kappa (\tau')\r \in A$.
The rest of the proof follows completely
those of [Gi, 4.18], where only the use of
[Gi, 4.13] there is replaced by the
following lemma:
\begin{lemma}
Let $\tau \le \delta$.  Suppose that $\l\l
A_i^{0\rho}, A^{1\rho}_i\r \mid \tau \le
\rho \le \delta, i < \kappa (\tau)^+ \r$
(where $\kappa (\tau) = \kappa^{+ \tau' +
n}$, if $\tau = \tau' + n$ for inaccessible
$\tau'$ and $n < \omega$, and $\kappa (\tau)
= \kappa^{+\tau + 1}$ otherwise) is an
increasing sequence of elements of ${\cal
P}'_{\ge \tau}$, satisfying the
following:\\
For every $i < \kappa (\tau)^+$ of
cofinality $\kappa (\tau)$
\begin{itemize}
\item[{\rm (a)}] $A_i^{0\tau} = \bigcup_{j<i}
A_j^{0\tau}$.
\item[{\rm (b)}] if $B \in A_i^{1\rho}$ then either $B
\supseteq A_i^{0\tau}$ or $B \in
\bigcup_{j<i} A_j^{1\rho}$.
\end{itemize}
\underline{Then} for every $i < \kappa (\tau)^+$
of cofinality $\kappa(\tau)$ and $B \in
\cup\{A_j^{1\rho}| j\ge i,\tau \le
\rho\le \delta\}$, either $B \supseteq
A_i^{0\tau}$ or there are
$\tilde i < i$, $\tau \le \tau_1 \le \tau_2
\le \cdots \le \tau_\ell \le \delta$ ($\ell
< \omega$) and $B_1 \in A_{\tilde
i}^{0\tau} \cap A_{\tilde
i}^{1\tau_1},\ldots, B_\ell \in A_{\tilde
i}^{0\tau} \cap A_{\tilde i}^{1\tau_\ell}$
such that for every $j, \tilde i\le j \le
i$
$$B \cap A_j^{0\tau} = B_1 \cap \cdots \cap
B_\ell \cap A_j^{0\tau}\ .$$

\end{lemma}

\medskip
\noindent
{\bf Proof.} Fix $i$ of cofinality $\kappa
(\tau)$. Let $B \in A_j^{1\rho}$ for some
$j$,~ $i\le j < \kappa (\tau)^+$ and $\rho$,
to $A_i^{0\tau} \in A_j^{0\tau}$, $\l\l
A_j^{0\nu}, A_j^{1\nu}\r \mid \tau \le \nu
\le \delta\r$ and $B \in A_j^{1\rho}$.
There are $\ell < \omega$, $\tau \le
\tau_1\le \cdots \le \tau_\ell \le \delta$
and $B_1 \in A_i^{0\tau} \cap
A_j^{1\tau_1},\ldots, B_\ell \in
A_i^{0\tau} \cap A_j^{1\tau_\ell}$ such
that
$$B \cap A_i^{0\tau} = B_1\cap \cdots \cap
B_\ell \cap A_i^{0\tau}\ .$$
By (a), there is $\tilde i < i$ such that
$B_1,\ldots, B_\ell \in A_{\tilde i}
^{0\tau}$.
Fix $k$, $1\le k \le \ell$. $B_k \in
A_{\tilde i}^{0\tau} \subseteq A_{\tilde
i}^{0\tau _k}$, so $B_k \subset A_{\tilde
i}^{0\tau_k}$ since they have the same
cardinality $\kappa (\tau_k) \subseteq
A_{\tilde i}^{0\tau_k}$.  Then by 1.2, $B_k
\in A_{\tilde i}^{1\tau_k}$.

Clearly, for every $i'$, $\tilde i \le
i'\le i$ $B \cap A_{i'}^{0\tau} = B \cap
(A_i^{0\tau} \cap A_{i'}^{0\tau}) = (B \cap
A_i^{0\tau}) \cap A_{i'}^{0\tau} =
((B_1\cap \cdots \cap B_k) \cap
A_i^{0\tau}) \cap A_{i'}^{0\tau} = (B_1
\cap \cdots \cap B_k) \cap (A_{i'}^{0\tau}
\cap A_i^{0\tau}) = B_1 \cap \cdots \cap
B_k \cap A_{i'}^{0\tau}$.\\
$\square$
\begin{lemma}
Let $\tau \le \delta$.Then the following
holds:
\begin{itemize}
\item [{\rm (a)}] if there is no inaccessible
$\tau'<\tau$ and $n < \omega$ such that
$\tau = \tau' + n$, then ${\cal P}_{<\tau}$
satisfies $\kappa^{+\tau +2}$-c.c. in
$V^{{\cal P}_{\ge \tau}}$
\item [{\rm (b)}] if $\tau = \tau' + n$ for some
inaccessible $\tau'< \tau$ and $0 < n <
\omega$, then ${\cal P}_{< \tau}$ satisfies
$\kappa^{+\tau + 1}$-c.c. in $V^{{\cal
P}_{\ge \tau}}$.
\item[{\rm (c)}] if $\tau$ is an inaccessible, then
${{\cal P}_{< \tau}}$ satisfies
$\kappa^{+\tau + 2}$-c.c. in $V^{{\cal
P}_{\ge \tau}}$.
\end{itemize}
\end{lemma}
The proof of this lemma repeats the proof
of 4.19 of [Gi].  We have here three cases
because of different cardinalities
according to the distance from an
inaccessible.

\begin{lemma}
The forcing ${\cal P}$ preserves all the
cardinals except probably the successors of
inaccessibles.
\end{lemma}
This follows from 1.11 and 1.12.
\begin{remark}
If one wants to preserve all the cardinals,
then instead of the full support taken
here, Easton type of support should be
taken.  Thus, fix some $\l\l
\underline{A}^{0\nu}, \underline{A}^{1\nu},
\underline F^\nu \r \mid \nu \le \delta \r
\in {\cal P}$.  Let $\underline{\cal P}$
consist of elements having Easton type support over the fixed
condition, i.e. $\l\l B^{0\nu}, B^{1\nu}$,
$G^\nu\r \mid \nu \le \delta \r$ will be in
$\underline{\cal P}$, iff for every
inaccessible $\lambda \le \delta$,\quad $|\{\nu|
\l B^{0\nu}, B^{1\nu}$, $G^\nu \r \ne \l
\underline A^{0\nu}, \underline A^{1\nu},
\underline F^\nu \r \}|< \lambda$.

Now we define our main forcing $\l {\cal
P}^{**}, \to \r$ as in [Gi]. Namely, let $G
\subseteq {\cal P}$ be generic. Set ${\cal
P}^{**} = \cup \{F^0 \mid \exists A^{00},
A^{10}, \l\l A^{0\tau}, A^{1\tau}, F^\tau\r
\mid 0 < \tau \le \delta\r$ $\l\l A^{0\nu},
A^{1\nu}, F^\nu\r \mid \nu \le \delta \r
\in G\}$. 
\end{remark}
The proof of the final lemma
repeats those of [Gi].

\begin{lemma}
In $V^{\cal P}$, $\l {\cal P}^{**}, \to \r$
satisfies $\kappa^{++}$-c.c.
\end{lemma}

\section{On gaps of size $\kappa^+$}
\setcounter{equation}{0}

The aim of the present section will be to
sketch the proof of the following:

\begin{theorem}
Let $\kappa$ be a cardinal of cofinality
$\omega$.  Suppose that for every $\nu <
\kappa$ the set $\{\alpha < \kappa\mid 
o(\alpha) \ge \alpha^{+ \nu}\}$ is cofinal
in $\kappa$.  Then there is a cofinality-preserving extension having the same
bounded subsets of $\kappa$ and  satisfying
$2^\kappa = \kappa ^{+\delta + 1}$ for
every $\delta < \kappa^{++}$.
\end{theorem}

For $\delta < \kappa$ it was done in [Gi],
so we concentrate on $\delta$'s between
$\kappa$ and $\kappa^{++}$.

Pick an increasing sequence $\kappa_0 <
\kappa_1 < \cdots < \kappa_n < \cdots <
\kappa$ so that
\begin{itemize}
\item[(a)] $\bigcup_{n<\omega} \kappa_n =
\kappa$
\item[(b)] $\kappa_n$ carries an extender
of the length $\kappa_n^{+ \kappa_{n-1}}$,
for every $n$, $0 < n < \omega$.
\end{itemize}

We force as in [Gi] but with the following
correspondence function $a_n$:
\begin{itemize}
\item[(i)] $\dom a_n \subseteq
\kappa^{+\kappa_{n-1}}$ and
\item[(ii)] $a_n(\kappa^{+\alpha + 1}) =
\kappa_n^{+\alpha + 1}$ for $\alpha \in
[\kappa^+, \kappa^{+\kappa_{n-1}}$).
\end{itemize}
This forcing produces an increasing (mod
finite) sequence of functions $\l f_\alpha
\mid \alpha \in [\kappa^+, \kappa^{+\kappa})\r$
such that $f_\alpha$ corresponds to $\kappa
^{+\alpha + 1}$.  Just define $f_\alpha(n)$
to be the element of the one element Prikry
sequence corresponding to
$\kappa_n^{+\alpha + 1}$, for every $n \ge
\min \{m\mid  \alpha \in [\kappa^{+m}, \kappa
^{+\kappa})\}$. The generic extension will
satisfy $2^\kappa \ge \kappa^{+\kappa}$.

In order to make $2^\kappa \ge
\kappa^{+\kappa^+}$ let us pick
$\kappa_n$'s carrying extenders of length $\kappa_n^{+\kappa_{n-1}^+}$. We require here only that $a_n
(\kappa^{+\alpha}) =
\kappa_n^{+a_{n-1}(\alpha)}$ and
$a_{n-1}(\alpha) \le \kappa_{n-1}^+$ for every $\alpha
\le \kappa^+$.  There is one minor problem
that for some $\alpha < \kappa^+$\enskip
$\kappa^{+\alpha + 1}$ can be in $\dom
f_n$, where the conditions at level $n$ are
of the form $\l a_n, A_n, f_n\r$.  We
required in [Gi] that $\dom a_n \cap \dom
f_n = \emptyset$.  Here we need to deal
with ordinals of cardinality
$\kappa^{+\alpha}$ and to add some of them
to $\dom a_n$.  We can either do it
directly and require $|a_n (\beta_1)| =
|a_n (\beta_2)|$ for any two such ordinals,
or explicitly add $\kappa^{+\alpha + 1}$ to
$\dom a_n$ and remove this value after a
nondirect extension is taken.  This
produces a generic extension satisfying
$2^\kappa \ge \kappa^{+\kappa^+}$.

Now in order to deal with arbitrary $\delta
\in [\kappa^+, \kappa^{++})$ we deal first
with $(\kappa^+)^m$ for every $m < \omega$
and then use the Rado-Milnor paradox, see
K. Kunen [Ku, Ch.~1, Ex.~20].

Fix $m < \omega$.  We use $\kappa_n$ which
carries an extender of  length
$\kappa_n^{+ (\kappa_{n-1}^+)^m}$.
We proceed as above with the following
addition:\\
for every $\kappa^+\le \beta <
(\kappa^+)^m$ let $\beta = (\kappa^+)^{k_1}
\cdot \alpha_1 + (\kappa^+)^{k_2}\cdot
\alpha_2 + \cdots + (\kappa^+)^{k_\ell} \cdot
\alpha_\ell$ where $\omega > k_1 > k_2 >
\cdots > k_\ell$, $\ell < \omega$,
$\alpha_1 ,\ldots \alpha_\ell < \kappa^+$,
then we require that $a_n (\beta) =
(\kappa_{n-1}^+ )^{k_1} \cdot a_n
(\alpha_1) + (\kappa_{n-1}^+)^{k_2}\cdot
a_n (\alpha_2) + \cdots +
(\kappa_{n-1}^+)^{k_\ell} \cdot a_n
(\alpha_\ell)$.

\section{Doing below the first inaccessible
above $\kappa$}
A straightforward application of the
techniques developed in [Gi] and here, can
get one above $\kappa^{+\kappa^{+\nu + 2}}$
for $\nu < \kappa$ starting with
$o(\kappa_n) = \kappa_n^{+\kappa_n^{+n +
\nu + 2}} + 1$ $(n < \omega)$, above
$\kappa^{+\kappa^{+\kappa^{+\nu + 2}}}$
starting with $o(\kappa_n) =
\kappa_n^{+\kappa_n^{\kappa_n^{+n + \nu +
2}}} + 1$,
etc.

Let us explain this dealing with
$\kappa^{+\kappa^{+\nu + 2}}$, i.e. we want
to have $2^\kappa \ge \kappa^{+\kappa^{+\nu +
2}}$ using an extender of  length
$\kappa_n^{+\kappa_n^{+n + \nu + 2}}$ for
each $n < \omega$.
The definition of the preparation forcing
${\cal P}$ is as 1.5 with the change in the
cardinals correspondence condition.  We
require the following:
\begin{itemize}
\item[(i)] $a_n (\kappa^{+ \tau + 1}) =
\kappa_n^{+n + \tau + 1}$ 
\item[(ii)] $a_n (\kappa^{+\kappa^{+\tau +
1}}) = \kappa_n^{+a_n (\kappa^{+\tau +
1})}$\newline
 for every $\tau \le \nu$ and $a_n$ as
in 1.5.
\end{itemize}
The rest of the construction is the same.
Notice only that in a previous section at
each $\kappa_n$ we had an extender of
inaccessible length which allowed us to have
many similar cardinals in the interval
$(\kappa_n$, the first inaccessible above
$\kappa_n$).  This in turn allowed us to pick
the correspondence between cardinals in the
intervals $(\kappa, \kappa^{+\delta})$ and
in $(\kappa_n$, the first inaccessible
above $\kappa_n)$ generically. In the
present situation the number of cardinals
to deal with is relatively small and so we
use (i) and (ii) above to define the
correspondence.  It remains [Gi, sec. 4,5],
where the correspondence was also
defined in advance.

Now we like to implement the Shelah idea
[Sh1] in order to show that below the first
inaccessible above $\kappa$ any gap is
possible, provided
that
for every $n<\ome$  we have an extender over
$\kap_n$ of the strength of the fixed point of
the aleph function above $\kap_n$  of order
$n$.  Let us first recall the definition. 

\begin{definition}
(Shelah [Sh1]).  Let
$C^0=$ the class of all infinite cardinals.
$C^{n+1}=\{\lam\in C^n\mid C^n\cap\lam$  has
order type $\lam\}$ and $C^\ome
=\bigcap_{n<\ome}C^n$.  The order of a cardinal
$\nu$  is the maximal $n\le \ome$ such that
$\nu\in C^n$. Elements of $C^n(n\ge 1)$ are
called fixed points of the $\aleph$-function of
order $n$.
\end{definition}

\begin{theorem}
 Suppose that
$\kap=\bigcup_{n<\ome}\kap_n$, $\kap_0 <\kap_1<\cdots
<\kap_n<\cdots ,o(\kap_n)=$ (the first fixed
point of the $\aleph$-function above $\kap_n$
of order $n$)$+1$ and there is no inaccessible
cardinal above $\kap$. Then for every $\lam$
there are cardinals and cofinalities preserving
the extension, not adding new bounded subsets to
$\kappa$ and  satisfying $2^\kap\ge\lam$.
\end{theorem}

\bigskip\noindent
{\bf Proof.} \quad Let $\mu >\kap$.  By S.
Shelah [Sh1], Lemma 2.5 there exists an
increasing sequence $\langle D_n\mid
n<\ome\rangle$ so that $\bigcup_{n<\ome}D_n=\{\chi|\chi$
is a cardinal $\kap^{++}\le\chi\le\mu^+\}$  and
for every $n<\ome$  there is no elements of
$C^n$  between $\kap$ and $\mu^+$ in a generic
extension $V_n$  of $V$  obtained by preserving
only elements of $D_n$ as cardinals between
$\kap^{++}$ and $\mu^+$.  Without loss of
generality, let us assume that $\kap^{++}$ and
$\mu^+$ are in $D_0$.  

We like to correspond elements of $D_n$ to the
cardinals between $\kap^{+n+2}_n$  and the
length of extender over $\kap_n$, i.e. the least
fixed point above $\kap_n$  of the $\aleph$-function
of order $n$.  Once we have this correspondence,
the rest of the construction will be as in the
previous section.

Let us define by induction on $n<\ome$ a two
place function $f_n$  from ordinals into
cardinals. For every ordinal $\xi$  set
$f_0(\xi,0)=\xi$ $f_0(\xi,1)=\xi^+$,
$f_0(\xi,\alp +1)=(f_0(\xi,\alp))^+$  and
$f_0(\xi,\alp)=\bigcup_{\bet <\alp}f_0(\xi,\bet)$
for a limit $\alp$. I.e.
$f_0(\xi,\alp)=\xi^{+\alp}$,  the $\alp$-th
cardinal past $\xi$.  Now define $f_1$.  First,
for an ordinal $\xi$  we define $f_1(\xi,0)$.
Set $\rho_0=\xi,\rho_1=f_0(\xi,\rho_0)$
and $\rho_{n+1}=f_0(\xi,\rho_n)$ for every
$n<\ome$.  Let
$f_1(\xi,0)=\bigcup_{n<\ome}\rho_n$. Then,
$f_1(\xi,0)$ will be the least fixed point of
the $\aleph$-function above $\xi$.
We define $f_1(\xi,1)$ in a similar fashion to
be the second fixed point of the $\aleph$-function
above $\xi$.  Thus, set $\rho_0=f_1(\xi,0),\rho_1=
f_0(f_1(\xi,0),\rho_0)$ and
$\rho_{n+1}=f_0(f_1(\xi,0),\rho_n)$  for every
$n<\ome$.  Set $f_1(\xi,1)=\bigcup_{n<\ome}\rho_n$.
Clearly, $f_1(\xi,1)=f_1(f_1(\xi,0),0)$.  Now
for every limit $\alp$ let $f_1(\xi,\alp)=\bigcup_{\bet
<\alp}f(\xi,\bet)$.  Define $f_1(\xi,\alp+1)$
to be least fixed point of the $\aleph$-function
above $f_1(\xi,\alp)$, i.e.
$f_1(\xi,\alp+1)=f_1(f_1(\xi,\alp),0)$.  

Suppose now that for every $k\le m$  $f_k$  is
defined and $f_k(\xi,\alp)$ is the $\alp$-th
fixed point of the order $k$ above $\xi$.
Define $f_{m+1}$.  First for any ordinal $\xi$
we define $f_{m+1}(\xi,0)$.  Set $\rho_0=f_m(\xi,0)$,   
$\rho_1=f_m(f_m(\xi,0),\rho_0)$  and
$\rho_{n+1}=f_m(f_m(\xi,0),\rho_n)$  for every
$n<\ome$.  Let $f_{m+1}(\xi,0)=\bigcup_{n<\ome}\rho_n$.
For a limit ordinal $\alp$  we set
$f_{m+1}(\xi,\alp)=\bigcup_{\bet
<\alp}f_m(\xi,\bet)$.  At successor stage let
$f_{m+1}(\xi,\alp +1)=f_{m+1}(f_{m+1}(\xi,\alp),0)$. 
This completes the inductive definition of $\langle
f_n\mid n<\ome\rangle$.

Now we are going to use such defined functions
$\langle f_n\mid n<\ome\rangle$ in order to produce
the desired correspondence between elements of
$D_n$'s and cardinals between $\kap^{+n+2}_n$
and the fixed point of order $n$  of the
$\aleph$-function above $\kap_n$. 

Fix $n<\ome$.  We may either work in the world
obtained by leaving only elements of $D_n$  to
be the cardinals in the interval $(\kap^{++},\mu^+)$
or relating the functions to $f_m(m<\ome)$  to the
set $D_n$ in the obvious fashion. Let us deal
with the first possibility.  Notation in this
case will be a bit simpler.  Thus, we assume
that the cardinals between $\kap^{++}$ and
$\mu^+$  are only the elements of $D_n$. Then,
by the choice of $D_n$, there is no fixed points
of the $\aleph$-function of order $n$  between
$\kap^{++}$ and $\mu^+$.  We define inductively
the correspondence $\pi$. 
For every $\alp\le\kap_n$  let $\pi(\alp)=\alp$.
For $\alp$, $\kap_n<\alp <\kap^+$ $\pi(\alp)$
may take any value in the interval
$(\kap_n,\kap^{+n+1}_n)$.  $\pi(\kap^+)=\kap^{+n+1}_n$
and $\pi(\kap^{++})=\kap^{+n+2}_n$.  Notice that
$\kap^{++}=f_0(\kap^{++},0)$. So, we have
defined the correspondence between sets $\{0,1\nek
\alp\nek \kap,\kap^+,\kap^{++}=f_0(\kap^{++},0)\}$
and $\{0,1\nek\alp\nek \kap_n,\kap^+_n\nek\kap^{+n}_n,
\kap^{+n+1}_n$, $\kap^{+n+2}_n=f_0(\kap^{+n+2}_n,0)$.
Using it we continue to the sets
$\{f_0(\kap^{++},0),f_0(\kap^{++},1)\nek
f_0(\kap^{++}\alp)\nek f_0(\kap^{++},f_0(\kap^{++},0))\}$
and $\{f_0(\kap^{+n+2}_n,0)$,
$f_0(\kap^{+n+2}_n,1)\nek f_0(\kap^{+n+2}_n,\alp)\nek
f_0(\kap^{+n+2}_n,f_0(\kap^{+n+2}_n,0))\}$.
Just $f_0(\kap^{++},\alp)$  will correspond to
$f_0(\kap^{+n+2}_n,\pi(\alp))$.  Set $\rho_1=
f_0(\kap^{++},f_0(\kap^{++},0))$ and
$\orho_1=f_0(\kap^{+n+2}_n,f_0(\kap^{+n+2}_n,0))$.  Then
$\pi(\rho_1)=\orho_1$.  We now consider the sets
$\{f_0(\rho_1,0),f_0(\rho_1,1)\nek
f_0(\rho_1,\alp)\nek f_0(\rho_1,\rho_1\}$ and 
$\{f_0(\orho_1,0),f_0(\orho_1,1)\nek
f_0(\orho_1,\alp)\nek f_0(\orho_1,\orho_1)\}$.
Extend $\pi$ to these sets by setting\break
$\pi(f_0(\rho_0,\alp))
=f_0(\orho_1,\pi(\alp))$.
Let $\rho_2=f_0(\rho_1,\rho_1)$ and $\pi(\rho_2)=\orho_2
=f_0(\orho_1,\orho_1)$.  We consider the sets
$\{f_0(\rho_2,0)\nek f_0(\rho_2,\rho_2)\}$ and
$\{ f_0(\orho_2,0)\nek f_0(\orho_2,\orho_2)\}$
Deal with them in the same fashion.
Continuing and using induction, be will be able
to extend $\pi$  up to $\rho_\ome =\bigcup_{n<\ome}
\rho_n$.  But, clearly,
$\rho_\ome=f_1(\kap^{++},0)$, i.e. the first
fixed point of the $\aleph$-function above
$\kap^{++}$. So we have the correspondence
between sets $\{0\nek\kap\nek
f_1(\kap^{++},0)\}$ and $\{0,\nek\kap_n\nek
f_1(\kap^{+n+2}_n,0)\}$.  Let us extend it to
the correspondence between the sets
$\{f_1(\kap^{++},0),f_1(\kap^{++},1)\nek f_1(\kap^{++},
\alp)\nek$ $f_1(\kap^{++}$, $f_1(\kap^{++},0)\}$  
and $\{f_1(\kap^{+n+2}_n,0),$\break
$f_1(\kap^{+n+2}_n,1)\nek
f_1(\kap^{+n+2}_n,\alp)\nek f_1(\kap^{+n+2}_n$,
$f_1(\kap^{+n+2}_n,0)\}$.  Just let $f_1(\kap^{++},\alp)$
correspond to $f_1(\kap^{+n+2}_n,\pi(\alp))$.  The
correspondence between the intervals $(f_1(\kap^{++},
\alp), f_1(\kap^{++},\alp+1))$ and
$(f_1(\kap^{+n+2}_n,\pi(\alp))$,  $f_1(\kap^{+n+2}_n,
\pi(\alp)+1))$  is defined using
$f_0(f_1(\kap^{++},\alp)$, as above.  Set
$\rho_1=f_1(\kap^{++},f_1(\kap^{++},0))$ and
$\orho_1=\pi(\rho_1)=f_1(\kap^{+n+2}_n,f_1
(\kap^{+n+2}_n,0))$.  We consider the sets $\{
f_1(\rho_1,0),f_1(\rho_1,1)\nek f_1(\rho_1,\alp)\nek
f_1(\rho_1,f_1(\rho_2,0))\}$ and $\{ f_1(\orho_1,0),
f_1
(\orho_1,1)
\nek f_1(\orho_1,\alp)$\break
$\nek
f_1(\orho_1,f_1(\orho_1,0))\}$.  Extend $\pi$ to
these sets as above.  Let $\rho_2=f_1(\rho_1,f_1
(\rho_1,0))$ and
$\orho_2=\pi(\rho_2)=f_1(\orho_1,f_1(\orho_1,0))$.
Continue in the same fashion by induction we
will obtain $\rho_{m+1}=f_1(\rho_m,f_1(\rho_m,0))$,
$\orho_{m+1}=\pi(\rho_{m+1})=f_1(\orho_m,f_1(\orho_m,0))$
and the extension of $\pi$  to the intervals
$[0,\rho_{m+1}],[0,\orho_{m+1}]$, where
$m<\ome$.  Set $\rho_\ome=\bigcup_{m<\ome}\rho_m$.
Then $\rho_\ome =f_2(\kap^{++},0)$ will be the
first fixed point of the $\aleph$-function of
the order $2$  above $\kap^{++}$.  By induction
it is easy to continue the process up to
$f_n(\kap^{++},0)$ (the first fixed point of the
$\aleph$-function of the order $n$  above
$\kap^{++}$.

This completes the inductive definition of
correspondence $\pi$.

The rest of the construction repeats those
of Section 1.

\end{document}